\documentclass[reqno]{amsart}
\usepackage{hyperref}
\usepackage[T1]{fontenc}

\def \Rd {{\mathbb{R}^{d}}}
\def \R {{\mathbb{R}}}

\def \R {{\mathbb{R}}}

\newcommand{\kk}{{\rm I~\hspace{-1.15ex}k}}
\newcommand{\II}{{\rm I~\hspace{-1.15ex}I}}
\newcommand{\jj}{{\rm I~\hspace{-1.15ex}k'}}

\newtheorem{Theo}{Theorem}[section]

\newtheorem{Lem}{Lemma}[section]

\newtheorem{Pro}{Proposition}[section]

\newtheorem{Rem}{Remark}[section]

\AtBeginDocument{
\vspace{8mm}}

\begin{document}
\title[\hfilneg   \hfil On a high-dimensional nonlinear spde]
{On a high-dimensional nonlinear stochastic partial differential equation}

\author[L. Boulanba,  M. Mellouk\hfil   \hfilneg]
{Lahcen Boulanba, Mohamed Mellouk}  

\address{Lahcen Boulanba\newline
 Universit\'{e} Cadi Ayyad, Facult\'{e} des Sciences
Semlalia, \\ D\'{e}partement de Math\'{e}matiques, B.P. 2390 Marrakech,
40.000, Maroc}
\email{l.boulanba@ucam.ac.ma}

\address{Mohamed Mellouk \newline
 MAP5, UMR CNRS 8145 \\
Universit\'{e} Paris-Descartes,45 rue des Saints-Pères
75270  Paris Cedex 06
France}
\email{mohamed.mellouk@parisdescartes.fr}

\thanks{Supported by l'action integrée MA/10/224 (volubilis)}
\keywords{Stochastic partial
differential equations; non-Lipschitz coefficients;
\hfill\break\indent
 martingale problem; Gaussian noises, Skorohod representation theorem}

\begin{abstract}
 In this paper we investigate a nonlinear stochastic partial differential equation (spde in short)
 perturbed by a space-correlated Gaussian noise in arbitrary dimension
 $d\geq1$, with a non-Lipschitz coefficient noisy term.
 The equation studied coincides in one dimension with the stochastic
 Burgers equation. Existence of a weak
 solution is established through an approximation procedure.
\end{abstract}

\maketitle
\numberwithin{equation}{section}
\newtheorem{theorem}{Theorem}[section]
\newtheorem{lemma}[theorem]{Lemma}
\newtheorem{proposition}[theorem]{Proposition}

\section{Introduction}

Let $d\geq 1$ be an
 integer and $D:\ = [0,1]^{d}$. In this paper, we study
 the following nonlinear stochastic partial differential equation on  $\mathbb{R}^{+}\times D$
\begin{eqnarray}\label{1}
 \frac{\partial u}{\partial t} (t,x)
  & = &  \frac{\partial}{\partial x_i}
  \left(\frac{\partial}{\partial x_i} u(t,x)+\frac{1}{2}u^{2}(t,x)\right)\nonumber  \\
  &&  + \,  \sigma(u(t,x)){\dot{F}}(t,x),\,\, t>0, \, x \in D, \, i=1,...,d
 \end{eqnarray}
 \noindent with Dirichlet boundary conditions
\begin{eqnarray}\label{2}
 u(t,x)= 0 ,\quad  t> 0 ,\quad x\in\partial D,
\end{eqnarray}
 \noindent and the initial condition
 \begin{eqnarray}\label{3}
 u(0,x)=u_{0}(x),\quad x\in D,
 \end{eqnarray}
\noindent where $u_{0}$\ is a continuous function on $D$ with values
in the interval $[0,1]$ and $\dot{F}$ is a noise on
$\mathbb{R}^{+}\times D$
  that is white in time and colored in space.
  The noisy coefficient term $\sigma$ is a real-valued
  function which will be described  precisely later. We will refer to the
  problem  $(\ref{1})-(\ref{3})$ by $Eq(d, u_0, \sigma)$. \\

\noindent  It is worth noting that
  the  equation $Eq(d, u_0, 0)$  gives the classic Burgers equation, which arises in physics and have extensively  been studied in the literature, see e.g.
   \cite{Ladyzhenskaya} and \cite{Matsumura and Nishihara}.\\

\noindent A class of equations of type $Eq(d, u_0, \sigma)$, with
noise depending only on time, was studied by Gy\"{o}ngy and Rovira
\cite{Gy Ro}. The authors took the space variable in a bounded
convex domain of $\mathbb{R}^{d}$, and the class investigated
contains  a version of the one-dimensional Burgers equation
\cite{BCJ, Bu, Da-De-Tem, dg, Gy, Gy-Nu, hopf} as a special case.
They proved
  existence, uniqueness and a comparison theorem under assumptions
  on coefficients which include global Lipschitz condition on the diffusion
   coefficient ~$\sigma$.\\

   \noindent
   However, the present work  is concerned with a non-Lipschitz diffusion
   coefficient. Since the $F$ is a space-time noise and
  $\sigma$ is non-Lipschitz the results of
  \cite{Gy Ro} cannot be applied here. \\

 \noindent The equation $Eq(1, u_0, \sigma)$  was studied by Kolkovska \cite{Kolk}.
 The author proved the existence of a weak solution by using an argument of tightness and solving a martingal problem.
  The case considered corresponds to the stepping stone model
 which describes migrating populations consisting of two types when the total
  population size does not change over time.
  For more details on
  the genealogy of a variation of this model, one can see
   Chapter 7 in \cite{Anja-these}.
 Our study extends this of \cite{Kolk} to
 a multidimensional case, and with a more general class of
 diffusion coefficients.

 \noindent It is also worth mentioning that a one dimensional
 version of (\ref{1}), perturbed with a space-time white noise, was investigated by Adler and Bonnet
 \cite{Adl-Bon, Bonnet} on the whole real line and with the
 non-Lipschitz noisy coefficient
 $\sigma(r) = \sqrt{r}$. The authors established the existence of a weak solution
 and its H\"older continuity. Moreover, they discussed (but did not prove)
 the uniqueness  of the solution.\\

\noindent In the present paper, we deal with a d-dimensional spdes. As it is well known, we can no longer use the space-time
white noise for the perturbation; indeed, in dimension $d\geq 2$, the fundamental solution associated with the
operator $\frac{\partial}{\partial
  t} - \Delta$ is not square integrable when $d\geq2$, while the one
    associated with the operator $\frac{\partial}{\partial
  t^2} - \Delta$ becomes less smooth as the dimension increases. Hence,  the martingale measure approach, introduced by Walsh
\cite{Walsh}, to investigate multidimensional SPDEs driven by a
space-time white noise gives solutions only in the space of random
distributions (see e.g. \cite{Albeverio-Haba-Russo,
 Dawson-Salehi,  Peszat-Zabczyk}). \\

\noindent In order to circumvent this difficulty, Dalang and Frangos
\cite{Dalang-Frangos} suggested to replace the space time white noise by some Gaussian noise which is
white in time and which has some space correlation in order to
obtain solutions in the space of real-valued stochastic processes.
See \cite{Dalang} and \cite{Marta-Sanz1}  and the references therein for more literature in
the subject.\\

\noindent The aim of this paper is to study the existence of
solution of a non-linear equation $Eq(d,u_0,\sigma)$ in arbitrary
dimension
 $d\geq1$, with a non-Lipschitz coefficient noisy term. The idea to prove our main result is to adopt an approach going back
  to Funaki \cite{Funaki}.

The paper is organized as follows. In the second section we give a
formulation of the problem. Section 3 is devoted to define a spatial
  discretization scheme of (\ref{1}) and obtain a system of
  stochastic differential equations. The existence and
  uniqueness of a unique strong solution for this system.
  Section 4 is devoted to the tightness of these approximating
  solutions. The existence of a weak solution of $Eq(d, u_0, \sigma)$
  is established in Section 5 by solving a martingale problem. We
  conclude with an Appendix containing some technical results needed
  along the paper.

  \noindent Note that all real positive constants are denoted by $c$
  regardless of their values and some of the standing parameters are not mentioned.

\section{Framework}
Let $(\Omega ,\mathcal{F},P)$ be  a complete probability space.  Let $ F = \{ F(\varphi),
\varphi\in{\mathcal{D}}(\mathbb{R}^{d+1}\} $ be a mean-zero $
\mathrm{L}^{2}(\Omega, \mathcal{F}, \mathbb{P})$-valued Gaussian
process with the covariance functional
\begin{eqnarray}\label{19}
\mathbb{E}(F(\varphi)F(\psi)) =
\int_{0}^{+\infty}ds\int_{D}dx\int_{D}dy\varphi(s,x)f(x,y)\psi(s,y),
\end{eqnarray}
\noindent where $f : \mathbb{R}^{d}\times
\mathbb{R}^{d}\longrightarrow \mathbb{R} $ denotes the correlation
kernel of the noise $F$.\\

\noindent Several authors discussed the weaker assumptions in such a way
that (\ref{19}) defines a covariance functional, see e.g. the
references \cite{Dalang-Frangos,
 {Dalang}, {Anja-these}}. \\

\noindent In this paper, we will take $f$ symmetric and
 belonging to $\mathcal{C}_b(D\times D)$, the space of continuous bounded functions. Following the same approach in \cite{Dalang}, the Gaussian process $F$
can be extended to a worthy martingale measure $M=\{M(t,A):=F([0,t]\times
A)\,:\,t\geq 0,\,A\in \mathcal{B}_{b}(D)\}$ which shall acts as
integrator, in the sense of Walsh \cite{Walsh}, where $\mathcal{B}_{b}(D)$ denotes the bounded Borel subsets of $D$. Let $\mathcal{%
G}_{t}$ be the completion of the $\sigma $--field generated by the random
variables $\{M(s,A),\;0\leq s\leq t,\;A\in \mathcal{B}_{b}(D)\}$%
. The properties of $F$ ensure that the process $M=\{M(t,A),\;t\geq 0,\;A\in
\mathcal{B}_{b}(D)\}$, is a martingale with respect to the
filtration $\{\mathcal{G}_{t}:t\geq 0\}$.
Then one can give a rigorous meaning  to solution of the formal equation
(\ref{1}).  A stochastic process
  $u : \Omega\times\mathbb{R}^{+}\times D\rightarrow\mathbb{R}$,
   which is jointly measurable and $\mathcal{%
G}_{t}$-adapted, is said to be
   a weak solution to the equation $Eq(d, u_0, \sigma)$ if  there exists a
   noise $F$ of the form (\ref{19})  such that
   for  each $\varphi\in\mathcal{C}^{2}(D)$ such that $\varphi = 0 $  on $\partial
   D$, and a.s. for almost all $t\geq 0$  and $x\in D$,
 \begin{eqnarray*}
 \int_{D}u(t,x)\varphi(x)dx
 & = & \int_D u_{0}(x)\varphi(x)dx
   +   \int_{0}^{t}\int_D u(s,x)\Delta\varphi(x)dxds\\
   & &
  -  \frac{1}{2}\int_{0}^{t}\int_D u^{2}(s,x)
       \sum_{i=1}^{d}\frac{\partial\varphi}{\partial x_{i}}(x)dxds
 + \int_{0}^{t}\int_D
 \sigma{(u(s,x))}\varphi(x)F(ds,dx).
 \end{eqnarray*}
   In order to formulate our main result, we assume the
 following hypothesis  :
\vspace{0.3cm}

\vspace{0.3cm} \noindent $\mathbf{(H)}$ \ \
  The function
$\sigma$
is H\"{o}lder continuous of order $ \frac{1}{2} \leq \alpha <
1$ i.e. there exists a positive constant $c$ such that for all $x, y
\in [0,1]$
  $$
  |\sigma(x) - \sigma(y)| \leq c |x-y|^{\alpha}.
  $$
Moreover, we assume that $\sigma$ is not identidically zero and
satisfies $\sigma(0) = \sigma(1) =0$.
 \begin{Theo}\label{18}
 Assume that $u_{0}$ is a continuous function on $D$ with values in
 the interval $[0,1]$, and that $\sigma$ satisfies $\mathbf{(H)}$. Then there
  exists a weak solution to the equation $Eq(d, u_0, \sigma)$.
 \end{Theo}
\begin{Rem}\label{rem2-ijmms-these}
 The hypothesis $\mathbf{(H)}$ implies the existence of
 a positive constant $c$ such that for all $x \in [0,1]$
  $$
  \sigma(x) \leq c   \min(x^{\alpha}, (1-x)^{\alpha}).
  $$
  \end{Rem}

\noindent\textbf{Examples}
\begin{description}
\item[1.]  The function $\sigma_1:  x \longmapsto \sqrt{x(1-x)}$, which corresponds
to the stepping stone model, satisfies $\mathbf{(H)}$ with $\alpha =
\frac{1}{2}$.

\item[2.]
Let $\gamma\in]\frac{1}{2}, 1]$. The function  $\sigma_2$ defined on
the interval $[0,1]$ by $\sigma_2(0) = 0$ \ and \ $\sigma_2(x) =
-x^{\gamma}\log(x)$ \ for \ $x\in]0,1]$  satisfies $\mathbf{(H)}$.
It is a concave function with values in $[0,1]$. Indeed, $\sigma_2$
is twice continuously differentiable on $]0,1[$ and $\sigma_2''(x)=
x^{\gamma - 2}(1 - 2\gamma + \gamma(1-\gamma)\log(x))$, for all
$x\in]0,1[$. Moreover, taking into account the continuity of
$\sigma_2$ on $[0,1]$ we deduce its concavity on $[0,1]$. It follows
that $|\sigma_2(x) - \sigma_2(y)| \leq \sup\{ |\sigma_2(x-y) -
\sigma_2(0)|, |\sigma_2(1) - \sigma_2(1-x+y)|\}$, for all $x, y \in
[0,1]$ such that $x\geq y$. Finally, we point out that
 $|\sigma_2(x-y) - \sigma_2(0)| = |\sigma_2(x-y)| \leq c
|x-y|^{\gamma-\epsilon}$, where $\epsilon \in ]0, \gamma
-\frac{1}{2}[$, \  and $|\sigma_2(1) - \sigma_2(1-x+y)| =
|\sigma_2(1-x+y)| \leq c |x-y|^{\eta}$ for any $\eta\in]0,1[$. The
result follows.

\end{description}

\section{The approximating processes}\label{sec3-3}
 \noindent Let $n\geq 2$ be an integer.  For $\kk:=(k_1, ..., k_d) \in \{1,..., n-1 \}^d$, \,  set $x_{\kk}^n=\left(\frac{k_{1}}{n},
\frac{k_{2}}{n},...,
  \frac{k_{d}}{n}\right)$ and define the interval
  $$\II^n_{\kk}=\displaystyle \Pi_{j=1}^d [\frac{k_j}{n}, \frac{k_j+1}{n}).$$
\noindent
 Consider the discretized
 version of $Eq(d, u_0, \sigma)$:
 \begin{align}\label{4}
 \left\{
 \begin{aligned}
& d u^{n}(t,x_{\kk}^n)
  =  \Delta^{n}u^{n}(t,x_{\kk}^n)dt
  + \frac{1}{2}
  \sum_{i=1}^{d}\nabla^{n}_{{i}}(u^{n}(t,x_{\kk}^n)^{2})dt
  +  n^{d}\sigma(u^{n}(t,x_{\kk}^n))F(dt,\II^n_{\kk}),\\
& u^{n}(0,x_{\kk}^n)\quad =  u_{0}({x_{\kk}^n}),
  \end{aligned}
  \right.
  \end{align}
where, the operator in the first term is given by
  \begin{eqnarray}
  \Delta^{n}u^{n}(t,x_{\kk}^n)
  &  = & n^{2}\sum_{i=1}^{d}\left(u^{n}(t,x_{\kk}^n+ \frac{e_{i}}{n}) +
  u^{n}(t,x_{\kk}^n- \frac{e_{i}}{n}) - 2u^{n}(t,x_{\kk}^n)\right),\nonumber \\ \noindent \mbox{and in second term by} \nonumber \\
   \nabla^{n}_{{i}}(u^{n}(t,x_{\kk}^n)^{2})
 & = & n\left(u^{n}(t,x_{\kk}^n+ \frac{e_{i}}{n})^{2}-
   u^{n}(t,x_{\kk}^n)^{2}\right),\nonumber
  \end{eqnarray}
\noindent with $\{e_i,  1\leq i\leq d \}$ denotes the canonical
basis of $\mathbb{R}^{d}$, and
$$
  F(t,\II^n_{\kk}) = \int_{0}^{t}\int_{\II^n_{\kk}}F(ds,dx).
  $$
The noises $\{F(t,\II^n_{\kk}),   t\geq 0\}$, derived from colored noise $F$ on  $ \mathbb{R}^{d}$ with covariance given by (\ref{19}),
  are one dimensional correlated Brownian
  motions. Moreover, the boundedness  of $f$ implies that
  $\langle F(\cdot, \II^n_{\kk}) \rangle_t \leq c  tn^{-2d}$,
  where $c$ is a constant depending on $f$. Indeed
  \begin{eqnarray*}
\langle F(\cdot, \II^n_{\kk}) \rangle_t
 & = &
\int_{0}^{t}\int_{\mathbb{R}^d}\int_{\mathbb{R}^d}
\mathbf{1}_{\II^n_{\kk}}(x)f(x,y)\mathbf{1}_{\II^n_{\kk}}(y)dsdxdy\\
& = & t \int_{\II^n_{\kk}}\int_{\II^n_{\kk}} f(x,y)dxdy\\
& \leq & c\, t\, n^{-2d},
  \end{eqnarray*}
where $  c = \max\{f(x,y); \ (x,y)\in D\times D\}$.\\

\noindent Now, we proceed as in \cite{MLuc} and we consider the
bijection between the grid $ D_{n}^{d}:=
\{x_{\kk}^n=\left(\frac{k_{1}}{n}, \frac{k_{2}}{n},...,
  \frac{k_{d}}{n}\right),\,\, (k_1, ..., k_d) \in \{1,..., n-1 \}^d\} $,
 which assigns to each $\left(\frac{k_{1}}{n}, \frac{k_{2}}{n}, ...
 \frac{k_{d}}{n}\right)$ the
 integer \ $k_{1} + (k_{2} - 1)(n-1) + ... + (k_{d} -
 1)(n-1)^{d-1}$.\newline

 \noindent Set $N := (n-1)^{d}$, and for each
 $i\in\left\{1,..,N\right\}$,\  let
 $\underline{i}^{n}$ be the unique element
 $\left(\frac{k_{1}}{n}, \frac{k_{2}}{n}, ...,\frac{k_{d}}{n}\right)$
  such that
  $i = k_{1} + (k_{2} - 1)(n-1) + ... + (k_{d} -
 1)(n-1)^{d-1}$ and $1\leq k_{j}\leq n-1$ for each $j\in\{1,..,d\}$.\\

\noindent  Finally, let
  $u^{n}_{i}(t)$ denotes  $u^{n}(t,\underline{i}^{n})$.
  Hence, the system (\ref{4}) can be written, for $i=1,...,d$,
\begin{align}\label{5}
 \left\{
 \begin{aligned}
& d u_{i}^{n}(t)  =  \sum_{j=1}^{N}[a_{i,j} u_{j}^{n}(t)
 +  \frac{1}{2}b_{i,j}u_{j}^{n}(t)^{2}(t)]dt +
 n^{d}\sigma(u_{i}^{n}(t))F(dt,\II_{\underline{i}^n}^{n})\\
 & u_{i}^{n}(0)  =  u_{0}(\underline{i}^{n}).
\end{aligned}
\right.
\end{align}
\noindent Let $A_{n}^{(d)}$  and $B_{n}^{(d)}$be the $N\times N $
matrices such that $A_{n}^{(d)} :\ = \left(a_{i,j}\right)_{1\leq i,
j \leq N}$ and $B_{n}^{(d)} :\ = \left(b_{i,j}\right)_{1\leq i, j
\leq N}$. \\

\noindent For $d=1$, it is well known that we have

\hspace{-5.5cm}
\begin{equation}
 a_{i,j} = \left\{
\begin{array}{rl}
n^{2}     & \mbox{if} \quad |i-j| = 1,      \\
-2n^{2}   & \mbox{if} \quad j=i, \\
 0        & \quad otherwise,
 \end{array}
 \right.
 \end{equation}
 \noindent and,

\hspace{-6cm}
 \begin{equation}\label{1spde-d}
  b_{i,j} = \left\{
\begin{array}{rl}
n     & \mbox{if} \quad j = i + 1,      \\
-n   & \mbox{if} \quad j=i, \\
 0        & \quad otherwise.
 \end{array}
 \right.
\end{equation}

 \noindent For $d \geq 2$, it is known from
\cite{MLuc} that $A_{n}^{(d)}$ can be obtained by induction using a
relation of recurrence and the fact that $A_{n}^{(1)}$ is known. For
the coefficients $b_{i,j}$, let us first introduce some notations. For an integer $\ell\geq 1$, let
$I_{\ell}$ be the $\ell\times \ell$ identity
matrix and set $M : = (n-1)^{d-1}$. One can easily show that  
%

$$
B_{n}^{(d)} = \left(
\begin{array}{cccccc}
B_{n}^{(d-1)} &   nI_{M}       &   0     &   \cdots   &  \cdots    &   0      \\
  0           &  B_{n}^{(d-1)} & nI_{M}   &    0       &  \cdots    &   \vdots \\
 \vdots       &   0            & \ddots   &  \ddots   &  \ddots    &   \vdots \\
 \vdots       & \cdots         & \ddots  &   \ddots   &  \ddots    &   0      \\
 \vdots       & \cdots         & \cdots  &   \ddots   &  \ddots    &  nI_{M}       \\
 0            & \cdots         & \cdots  &   \cdots   &   0        &  B_{n}^{(d-1)}
 \end{array}
 \right) - nI_{N}.
$$
\noindent Note that  $B_{n}^{(1)}$ is the $(n-1)\times (n-1)$-matrix
defined by (\ref{1spde-d}).\\
\noindent The main result of this section is the
following
\begin{Pro}\label{pro4-1}
 For any integer $n\geq 1$ and for any initial random condition $u^{n}(0)=
 (\xi_{1}^{n},...,\xi_{N}^{n})\in[0,1]^{N}$, the system
\begin{align}\label{21}
\left\{
\begin{aligned}
 & d u_{i}^{n}(t)  =  \sum_{j=1}^{N}[a_{i,j} u_{j}^{n}(t)
 +  \frac{1}{2} b_{i,j}u_{j}^{n}(t)^{2}]dt +
 n^{d}\sigma{(u_{i}^{n}(t))}F(dt,\II_{\underline{i}^{n}}^{n}), \\
& u_{i}^{n}(0)\  =  \xi_{i}^{n},
\end{aligned}
 \right.
\end{align}
 admits a unique strong solution
 $u^{n}(.)\in \mathcal{C}([0,+\infty[,[0,1]^{N})$.
\end{Pro}
\noindent {\textbf{Proof}}.

We consider first the following system
\begin{align}\label{22}
\left\{
\begin{aligned}
& d u_{i}^{n}(t)  =   \sum_{j=1}^{N}\left[a_{i,j} u_{j}^{n}(t)
 + \frac{1}{2}b_{i,j}g(u_{j}^{n}(t))\right] dt +
 h(u_{i}^{n}(t))F(dt,I_{i}^{n}),\\%
&  u_{i}^{n}(0)  =  \xi_{i}^{n},
\end{aligned}
 \right.
\end{align}
 \noindent where $h : \mathbb{R}\longrightarrow\mathbb{R}$ is
defined by $h(x) = n^{d}\sigma(x)\mathbf1_{\{ 0 \leq x
\leq 1 \}}$, and $g : \mathbb{R}\longrightarrow\mathbb{R}$ is
defined by $g(x) = x^{2}\mathbf1_{\{-1 \leq x \leq 1
\}}$. In order to show that the continuous coefficients of the
finite dimensional system (\ref{22}) satisfy the growth linear
condition, we write it under
 the following form
\begin{align}\label{22-cor}
\left\{
\begin{aligned}
& d U^{n}(t)  =  A_{n}^{(d)}U^{n}(t)dt +
\frac{1}{2}B_{n}^{(d)}\Sigma(U^{n}(t))dt +
 K(U^{n}(t))dF_t\\%
&  U^{n}(0)  =  u^{n}(0),
\end{aligned}
 \right.
 \end{align}
with \ $\displaystyle  U^{n}(t): = \left(u_i^n(t),   1\leq i\leq
N\right), \ \
 \Sigma(U^{n}(t)): = \left(g(u_i^n(t)),   1\leq i\leq N\right),$ \
$\displaystyle  K(U^{n}(t)): = diag(h(u_i^n(t)), 1\leq i\leq N), \
\textit{and}\ \  F_{t}: = \left(F(t, \II_{\underline{i}^{n}}^{n}), \ 1\leq i\leq
N\right)$.\\

\noindent Since the first term in the right hand side of
(\ref{22-cor}) is linear, then it satisfies the linear growth
condition. For the second one, using an elementary property of
matrix norms and the equivalence of norms on $\mathbb{R}^{N}$, we
can write
\begin{eqnarray*}
\| B_{n}^{(d)}\Sigma (U^{n}(t))\| & \leq &
\|B_{n}^{(d)}\|\|L(U^{n}(t))\|\\
& \leq  & c \sum_{i=1}^{N}|u_i^n(t)|
  \leq   c \|U^{n}(t)\|,
\end{eqnarray*}
where the constant $c$ depends on $n$ and $d$. Concerning the noise
coefficient, owing to Remark \ref{rem2-ijmms-these} we obtain  a.s. for $t\geq 0$,
\begin{eqnarray*}
\sum_{i=1}^{N}h(u_i^n(t))^{2} & \leq & c n^{2d}\sum_{i=1}^{N}
u_i^n(t)^{2\alpha}\mathbf1_{\{ 0 \leq u_{i}^{n}(t) \leq 1
\}}\\
& \leq & c n^{2d}\sum_{i=1}^{N} (1 + (u_i^n(t))^{2})\\
& \leq & c n^{2d} N \left(1 + \|U^n(t)\|^{2}\right).
\end{eqnarray*}
Thus, the existence of a weak solution with continuous sample paths
for  (\ref{22})  is given by Theorem \ref{th4-app} in the Appendix.\\

\noindent Now, we will show that $0 \leq u_{i}^{n}(t)\leq 1 $ for
all $i = 1,...,n$ and $t\geq 0$.
 Taking $\rho(x) = x^{\beta}$ with $1 \leq \beta \leq
2\alpha$, and using the hypothesis $\mathbf{(H)}$ and the continuity
of the trajectories of $u_{i}^{n}$, we get fot $0\le t \le T$
\begin{eqnarray*}
\int_{0}^{t}\frac{\mathbf1_{\{u_{i}^{n}(s) > 0
\}}}{\rho(u_{i}^{n}(s))}d\langle u_{i}^{n} \rangle_{s} & = & n^{2d}
\int_{0}^{t}\frac{\mathbf1_{\{u_{i}^{n}(s) > 0
\}}h^{2}(u_{i}^{n}(s))}{\rho(u_{i}^{n}(s))}ds\\
& \leq & n^{2d}T \\
& < &  +\infty.
\end{eqnarray*}
Then, by Lemma \ref{3-7} we obtain that the local time,
$L_{t}^{0}(u_{i}^{n})$, of the semi-martingale $u_{i}^{n}$ satisfies
$L_{t}^{0}(x_{i}^{n}) =0$
for all $t>0$, a.s, and $i = 1, \cdots, N$.  \\

\noindent Moreover, set $ x^{-} :\ = \max(0,-x)$ and by applying the
Tanaka formula and summing over all indices $i=1,..,N$, we get for $t\geq 0$,
\begin{eqnarray*}
 \sum_{i =1}^{N}u_{i}^{n}(t)^-
&  = & -\int_{0}^{t}\sum_{i = 1 }^{N}\mathbf1_{\{u_{i}^{n}(s)\leq
0\}} \sum_{j = 1}^{N}a_{i,j}u_{j}^{n}(s)ds -
\frac{1}{2}\int_{0}^{t}\sum_{i = 1}^{N}\mathbf1_{\{u_{i}^{n}(s) \leq
0\}} \sum_{j = 1}^{N}b_{i,j}g(u_{j}^{n}(s))ds\\
& :=&  I_{1}  + I_{2}.
\end{eqnarray*}

\noindent
Since  $-x \leq x^{-}$ and taking into account the fact that $a_{i,j} \geq 0$ for all $j\neq
i$ and $ \sum_{ i = 1, i\neq j}^{N}a_{i,j} \leq 2 d n^{2} $

\begin{eqnarray}
\label{Iun}
I_{1} &=&-\int_{0}^{t}\sum_{j = 1}^{N}u_{j}^{n}(s)\sum_{ i =
1}^{N}a_{i,j}\mathbf1_{\{u_{i}^{n}(s)\leq 0\}}ds \nonumber\\
&=&-\int_{0}^{t}\sum_{j = 1}^{N}\big[u_{j}^{n}(s) \big(-2 n^2 \mathbf1_{\{u_{j}^{n}(s)\leq 0\}} +
\sum_{ i = 1, i\neq
j}^{N}a_{i,j}\mathbf1_{\{u_{i}^{n}(s)\leq 0\}}\big)\big]\nonumber\\
&=&-\int_{0}^{t}\sum_{j = 1}^{N}\big[u_{j}^{n}(s)^{-}\, 2 n^2 \mathbf1_{\{u_{j}^{n}(s)\leq 0\}} +
u_{j}^{n}(s) \sum_{ i = 1, i\neq
j}^{N}a_{i,j}\mathbf1_{\{u_{i}^{n}(s)\leq 0\}}\big]\nonumber\\
&=&\int_{0}^{t}\sum_{j = 1}^{N}\big[u_{j}^{n}(s)^{-}\, (- 2 n^2) \mathbf1_{\{u_{j}^{n}(s)\leq 0\}} -
u_{j}^{n}(s) \sum_{ i = 1, i\neq
j}^{N}a_{i,j}\mathbf1_{\{u_{i}^{n}(s)\leq 0\}}\big]\nonumber\\
&\le&\int_{0}^{t}\sum_{j = 1}^{N}u_{j}^{n}(s)^{-}\left(\sum_{ i =
1}^{N}a_{i,j}\mathbf1_{\{u_{i}^{n}(s)\leq 0\}} \right) ds \le 2 d
n^{2} \int_{0}^{t}\sum_{i=1}^{N} u_{i}^{n}(s)^{-} ds.
\end{eqnarray}

\noindent
 For $I_{2}$, notice first that by using the definition of $g$,
\begin{eqnarray*}
I_{2} = \frac{1}{2}\int_{0}^{t}\sum_{i=1}^{N}\mathbf1_{\{-1\leq
u_{i}^{n}(s) \leq 0\}} \left(nd u_{i}^{n}(s)^{2} - \sum_{j=1, j\neq
i}^{N}b_{i,j}u_{j}^{n}(s)^{2} \right)ds.
\end{eqnarray*}
Taking into account the fact that $b_{i,j} \geq 0$ for all $j\neq
i$, we infer that
\begin{eqnarray*}
I_{2} \leq \frac{1}{2}nd\int_{0}^{t}\sum_{i=1}^{N}\mathbf1_{\{-1\leq
u_{i}^{n}(s) \leq 0\}} u_{i}^{n}(s)^{2} ds.
\end{eqnarray*}
Moreover,  for  $0 \leq s\leq t$ such that $-1 \leq u_{i}^{n}(s)
\leq 0$, we have
\begin{eqnarray*}
u_{i}^{n}(s)^2 = (- u_{i}^{n}(s))^2 \leq - u_{i}^{n}(s) =
u_{i}^{n}(s)^{-}.
\end{eqnarray*}
Hence,
\begin{eqnarray}\label{24-cor}
I_{2} \leq \frac{1}{2}nd\int_{0}^{t}\sum_{i=1}^{N} u_{i}^{n}(s)^{-}
ds.
\end{eqnarray}
From the estimates (\ref{Iun}) and (\ref{24-cor}), we deduce
\begin{eqnarray*}
\sum_{i =1}^{N}u_{i}^{n}(t)^- \leq
c\,\int_{0}^{t}\sum_{i=1}^{N} u_{i}^{n}(s)^{-} ds.
\end{eqnarray*}
 Then, by Gronwall's Lemma, we obtain
that $\sum_{i = 1}^{N} u_{i}^{n}(t)^{-} = 0 $. Consequently, the
solution is non-negative for each $t\geq 0$.\\
\noindent
 Concerning the fact that  $u_i^n(t) \leq 1$, it can be
easily checked by using the same arguments as above for  $(1 -
u_i^n(t))$.\\

\noindent Therefore, the system (\ref{21}) admits a weak solution
with trajectories in $\mathcal{C}([0,+\infty[,[0,1]^{N})$.\\

\noindent To show the pathwise uniqueness for (\ref{21}), we assume
that $u^{1,n}(t): \  = (u_i^{1,n}(t), 1\leq i\leq N )$ and
$u^{2,n}(t): \ = (u_i^{2,n}(t), 1\leq i\leq N )$ are two weak
solutions to (\ref{21}), with the same noise $F$ and the same
initial data.\\

\noindent Set $y_{i}^{n}: =  u_i^{1,n} - u_i^{2,n}$, then for $t>0$,
\begin{eqnarray*}
y_{i}^{n}(t) & = & \sum_{j = 1}^{N}a_{i,j}\int_{0}^{t} y_{j}^{n}(s)
ds + \frac{1}{2}\sum_{j = 1}^{N}b_{i,j}\int_{0}^{t}\left(
u_j^{1,n}(s)^{2} - u_j^{2,n}(s)^{2}\right)ds\\
& & + n^{d}\int_{0}^{t}\left[\sigma(u_{i}^{1,n}(s)) -
\sigma(u_{i}^{2,n}(s))\right] F(ds,\II_{\underline{i}^{n}}^{n}).
\end{eqnarray*}
Consequently, by the boundedness of the correlation kernel $f$,
\begin{eqnarray*}
\langle y_{i}^{n}(\cdot) \rangle_{t} \leq  c
\int_{0}^{t}\left[\sigma(u_{i}^{1,n}(s)) -
\sigma(u_{i}^{2,n}(s))\right]^{2}ds.
\end{eqnarray*}
Owing to Remark \ref{rem2-ijmms-these} and Lemma \ref{3-7}, we
obtain that a.s.
 $\displaystyle  L_{t}^{0}\left(y_{i}^{n}\right) = 0$, for all
 $i\in\{1,...,N\}$. Applying the Tanaka's formula for the
 continuous semimartingale $y_{i}^{n}$, it follows that
 \begin{eqnarray*}
|y_{i}^{n}(t)| & = &  \sum_{j  =
1}^{N}a_{i,j}\int_{0}^{t}sgn(y_{j}^{n}(s))y_{j}^{n}(s)ds\\
& &  + \frac{1}{2}\sum_{j =
1}^{N}b_{i,j}\int_{0}^{t}sgn(y_{j}^{n}(s))\left(
u_j^{1,n}(s)^{2} - u_j^{2,n}(s)^{2}\right)ds\\
& & + n^{d}\int_{0}^{t}sgn(y_{i}^{n}(s))\left[\sigma(u_{i}^{1,n}(s))
- \sigma(u_{i}^{2,n}(s))\right] F(ds,\II_{\underline{i}^{n}}^{n}).
 \end{eqnarray*}
 Noting that the sign function is defined by
 \begin{equation*}
 sgn(x) = \left\{
\begin{array}{l}
+1, \quad x\geq 0;\\
-1, \quad x < 0.
\end{array}
\right.
\end{equation*}
By summing over all $i\in\{1,...,N\}$ and taking the expectation, we
get
\begin{eqnarray}
\mathbb{E} \left(\sum_{i =1}^{N}|y_{i}^{n}(t)|\right) & = & \mathbb{E}\left( \sum_{i, j
=
1}^{N}a_{i,j}\int_{0}^{t}sgn(y_{j}^{n}(s))y_{j}^{n}(s)ds\right)\nonumber\\
& & + \frac{1}{2}\mathbb{E} \left(\sum_{i, j =
1}^{N}b_{i,j}\int_{0}^{t}sgn(y_{j}^{n}(s))\left(u_j^{1,n}(s)^{2}
- u_j^{2,n}(s)^{2}\right)ds\right)\nonumber\\
& :=&  J_{1} + J_{2}. \label{5-uni-1}
\end{eqnarray}
For $J_{1}$, the boundedness of the coefficients $a_{i,j}$ yields
\begin{eqnarray}
J_{1} & = & \mathbb{E} \int_{0}^{t}\sum_{ j =
1}^{N}sgn(y_{j}^{n}(s))y_{j}^{n}(s) \left(\sum_{ i
= 1}^{N}a_{i,j}\right)  ds \nonumber\\
& \leq & c \,\mathbb{E}\int_{0}^{t}\sum_{ j = 1}^{N}|y_{j}^{n}(s)|ds,
\label{5-uni-2}
\end{eqnarray}
where $c$ is a constant depending on $n$. Concerning $J_{2}$, using
the boundedness of the coefficients $b_{i,j}$ and the fact that the
solutions are with values in the interval $[0,1]$, we can write
\begin{eqnarray}
J_{2}& = &\frac{1}{2}\mathbb{E}\int_{0}^{t}\sum_{ j =
1}^{N}sgn(y_{j}^{n}(s)) y_{j}^{n}(s)(u_{i}^{1,n}(s) +
u_{i}^{2,n}(s))  \left(\sum_{ i
= 1}^{N}b_{i,j}\right)ds \nonumber\\
& \leq & c \mathbb{E}\int_{0}^{t}\sum_{j = 1}^{N}|y_{j}^{n}(s)|ds.
\label{5-uni-3}
\end{eqnarray}
Combining the estimates (\ref{5-uni-1})--(\ref{5-uni-3}) and
applying Gronwall's Lemma, we obtain for all $t\geq 0$
$$
\mathbb{E}\sum_{i=1}^{N}\left|y_i^{n}(t)\right| = 0.
$$

\noindent Finally, the well-known theorem of Yamada and Watanabe
\cite{Yam-Wat} (see also Corollary 3.23 in \cite{Karatzas}) implies
that (\ref{21})
 has a unique strong solution. \hfill $\square$
\section{The tightness of the approximating processes}
\noindent As in \cite{Sturm 1}, let $Y^{n}$ be a simple random walk
 whose generator is the discrete
Laplacian $\Delta^n$ on the lattice $D_{n}^{d}$. To get the system (\ref{5}) in its mild form we define the
fundamental solution $p^n_d $ associated with $\Delta^n$. \\

\noindent For $
\underline{i}^n, \underline{j}^n \in D_{n}^{d}$,
set
$$
p^n_d(t, \underline{i}^n, \underline{j}^n) = n^d\mathbb{P}(Y^{n}_{t}
= \underline{j}^n\mid  Y^{n}_{0} = \underline{i}^n).
$$
 The system (\ref{5}) can be now written in its variation of
constant form
\begin{eqnarray*}
u_i^n(t) & = & \frac{1}{n^d} \sum_{j =1}^{N} p^n_d(t,
\underline{i}^n, \underline{j}^n)u_0(\underline{j}^n) +
\frac{1}{2n^d}\int_0^t \sum_{j =1}^{N} p^n_d(t
-s, \underline{i}^n, \underline{j}^n)b_{i,j}u_j^n(s)^2ds\\
      &  & +  \,\int_0^t \sum_{j =1}^{N} p^n_d(t
-s, \underline{i}^n, \underline{j}^n)\sigma(u_j^n(s))F(ds, \II_{\underline{i}^{n}}^{n}).
\end{eqnarray*}
\noindent In the remainder of this paper, we denote $u_i^n(t)$ by
$u^n(t,\underline{i}^n)$, where the relationship between $i$ and
 $\underline{i}^n$ is described in the previous section. Moreover,
 in order to define $u^n(t,\cdot)$ for all $x\in
D$, we set $\kappa_n(x) = \underline{i}^n$ for $x\in I_i^n$, and
$\bar{p}^n_d(t , x, y) =  p^n_d(t , \kappa_n(x), \kappa_n(y))$ for
$x, y \in D$.
 Using a polygonal interpolation like as in
\cite{Sturm 1}, we can write for all $x\in D$
\begin{eqnarray}
u^n(t,x) & = &  \int_D p^n_d(t, x, \kappa_n(y))u_0(\kappa_n(y))dy +
\frac{1}{2} \int_0^t
\int_D p^n_d(t-s , x, \kappa_n(y))u^n(s,y)^2dyds\nonumber\\
         &  & + \int_0^t \int_D
p^n_d(t-s , x, \kappa_n(y))\sigma(u^n(s,y))F(ds, dy)\nonumber\\
& \equiv &  u^n_1(t,x) + u^n_2(t,x) + u^n_3(t,x).\label{24}
\end{eqnarray}
\noindent The main result of this section is as follows.
\begin{Pro} For each $T>0$, the sequence $\{u^n, n\in \mathbb{N}\}$ is tight in
the space $\mathcal{C}\left([0,T], \mathcal{H}\right)$, where
$\mathcal{H} :\ = \mathcal{C}\left(D, [0,1]\right)$.
\end{Pro}
\noindent \textbf{Proof.} Note firstly that
$u^n(t,x)$ are in $[0,1]$. Consequently, there exists a constant
$c(p,T)>0$ such that for any $p\geq 1$,
\begin{eqnarray}\label{10}
\sup_n \mathbb{E}\left[\sup_{0\leq t \leq T}\sup_{x\in
D}|u^n(t,x)|^p\right] \leq c(p,T).
\end{eqnarray}
\noindent To handle the stochastic integral, we will use the
factorization method introduced by Da Prato el al. in \cite{Daprato}. Indeed, it is well known that one can write
for $\beta\in (0,1)$,
\begin{eqnarray}\label{12}
&& \int_0^t \int_D \bar{p}^n_d(t-s , x, y)\sigma(u^n(s,y))F(ds,
dy)\nonumber\\
&& = \frac{\sin(\beta \pi)}{\pi} \int_0^t \int_D (t-s)^{\beta
-1}\bar{p}^n_d(t-s , x, y)Y(s,y)dyds,
\end{eqnarray}
 where
$$
 Y(s,y) \ : \ = \int_0^s\int_D
(s-r)^{-\beta}\bar{p}^n_d(s-r , y, z)\sigma(u^n(r,z))F(dr,dz).$$ In
the sequel, we will take $0<\beta<\frac{1}{2}$. Hence, by
Burkholder's inequality, the boundedness of $\sigma$ and $f$ and the
estimate $i)$ of Lemma \ref{9}, we obtain for any $T>0$ and $p>1$,
that there exists a constant $c(T,p)
> 0$ such that
\begin{eqnarray}\label{11}
\sup_{s\in [0,T]}\sup_{y\in D} \mathbb{E}|Y(s,y)|^{p} \leq c(p,T).
\end{eqnarray}
\noindent Let us now prove an estimate on the increments in both the
space and time variables. Firstly, we are concerned with the third
term in (\ref{24}). \noindent Clearly, for all $x,  y \in D$ and
$0\leq s\leq t \leq T$
\begin{eqnarray*}\label{13}
\mathbb{E}|u^n_3(t,x) - u^n_3(s,y)|^{p} \leq
c\left(\mathbb{E}|u^n_3(t,x) - u^n_3(t,y)|^{p} +
\mathbb{E}|u^n_3(t,y) - u^n_3(s,y)|^{p}\right).
\end{eqnarray*}
\noindent By Burkholder's inequality and the fact that $\sigma$ and
the kernel noise $f$ are bounded, we can write
\begin{eqnarray}\label{17}
\mathbb{E}|u_3^n(t,x) - u_3^n(t,y)|^{p} & = &
\mathbb{E}\left|\int_{0}^{t}\int_{D}\left[\bar{p}^n_d(t -s,x,z) -
\bar{p}^n_d(t - s, y,
z)\right]\sigma(u^n(s,z))F(dz,ds)\right|^{p}\nonumber\\
& \leq & c \left(\int_{0}^{t}\left(\int_{D}\left|\bar{p}^n_d(t
-s,x,z) - \bar{p}^n_d(t - s, y,
z)\right|dz\right)^{2}ds\right)^{p/2}\nonumber\\
& \leq & c \|x-y\|^{\left(\frac{1}{q} -
\frac{1}{2}\right)\frac{p}{2}},
\end{eqnarray}
\noindent where $ 1 < q < 2$. Note that in the last inequality, we
have used the estimate $ii)$ of Lemma ~~\ref{9}. On the other hand,
we have
\begin{eqnarray*}
&& \mathbb{E}|u^n_3(t,y) - u^n_3(s,y)|^{p}\\
&& = \mathbb{E}\left|\int_0^t\int_D\left[(t-r)^{\beta
-1}\bar{p}^n_d(t-r,y,z) -
 (s-r)^{\beta-1}\bar{p}^n_d(s-r,y,z)\mathbf1_{[0,s]}(r)\right]Y(r,z)dzdr\right|^{p}\\
 && \leq
c\left(\mathbb{E}\left|\int_0^s\int_D(s-r)^{\beta
-1}(\bar{p}^n_d(t-r,y,z)-\bar{p}^n_d(s-r,y,z))Y(r,z)dzdr\right|^{p}\right.\\
&& \quad + \quad \mathbb{E}\left|\int_0^s\int_D((t-r)^{\beta -1} -
(s-r)^{\beta -1})
\bar{p}^n_d(t-r,y,z)Y(r,z)dzdr\right|^{p}\\
&&\quad  + \quad \left.\mathbb{E}\left|\int_s^t\int_D(t-r)^{\beta
-1}
\bar{p}^n_d(t-r,y,z)Y(r,z)dzdr\right|^{p}\right)\\
&& \equiv  A_1 +  A_2 + A_3.
\end{eqnarray*}
\noindent By H\"{o}lder's inequality, (\ref{11}) and
 $i)$ of the Lemma \ref{9}, we get
\begin{equation}
 A_3 \leq c |t-s|^{(\beta - \lambda d)p},
\end{equation}
\noindent  where $0 < \beta <\frac{1}{2} $  and $0 < \lambda
<\frac{\beta}{d}$. Concerning $A_2$, the same arguments as above
yield
\begin{eqnarray}
  A_2 & \leq & c\left(\int_{0}^{s}[(s-r)^{\beta -1} -
              (t-r)^{\beta -1}](t-r)^{-\lambda d}dr \right)^{p}\nonumber\\
      & \leq & c |t-s|^{-2\lambda d p} \left(\int_{0}^{s}[(s-r)^{\beta-1}
               - (t-r)^{\beta-1}]dr\right)^{p} \nonumber\\
      & \leq & c |t-s|^{p(\beta-\lambda  d )}.
\end{eqnarray}
\noindent For $A_1$, we use firstly Burkholder's inequality
\begin{eqnarray*}
A_1 & = &
\mathbb{E}\left|\int_{0}^{s}\int_{D}(\bar{p}^n_d(t-r,y,z)-\bar{p}^n_d(s-r,y,z))
          \sigma(u^n(r,z))F(dr,dz)\right|^p\\
    & \leq & c \mathbb{E} \left( \int_{0}^{s}\int_{D}\int_{D}(\bar{p}^n_d(t-r,y,z)-\bar{p}^n_d(s-r,y,z))
          \sigma(u^n(r,z))f(z,x)\right.\\
  & &    \left. \times (\bar{p}^n_d(t-r,y,x)-\bar{p}^n_d(s-r,y,x))
       \sigma(u^n(r,x))drdzdx\vphantom{L\int_{0}^{s}}\right)^{p/2},\label{23}
\end{eqnarray*}
\noindent and since $\sigma$ and $f$ are bounded, and by $iv)$ of
the Lemma \ref{9}
\begin{eqnarray}
A_1 & \leq & c\left( \int_{0}^{s}\int_{D}\int_{D}
             (\bar{p}^n_d(t-r,y,z)-\bar{p}^n_d(s-r,y,z))|z-x|^{-\alpha}\right.\\
       & & \times\left.
              (\bar{p}^n_d(t-r,y,x)-\bar{p}^n_d(s-r,y,x))\right)^{p/2}\nonumber\\
    & \leq & c \left(\int_{0}^{s}\| \bar{p}^n_d(t-r,y,\cdot)-\bar{p}^n_d(s-r,y,\cdot)\|
           _{\alpha}^{2}dr\right)^{p/2}\nonumber\\
    & \leq & c |t-s|^{(1-\frac{\alpha}{2})\frac{p}{2}}.
\end{eqnarray}
\noindent Therefore, taking into account (\ref{17})-(\ref{23}),
there exists a positive constant $c$  and $0 < \delta_1<
\frac{1}{2}$ and $0 < \delta_2 < \frac{1}{4}$ such that
\begin{eqnarray}\label{15}
\mathbb{E}|u^n_3(t,y) - u^n_3(s,y)|^{2p} \leq c
\left(|t-s|^{\delta_1 p} + \|x-y\|^{\delta_2 p}\right).
\end{eqnarray}
 From the Proposition 3.2 in \cite{MLuc}, we deduce that (\ref{15})
holds for $u_1^n$. Concerning $u_2^n$, using H\"{o}lder's
inequality, (\ref{10}) and the estimates on $\bar{p}^n_d$ given by
the Lemma \ref{9}, we deduce that (\ref{15}) holds for $u_2^n$.
Finally, we deduce the tightness of the sequence $\{ u^n, n\geq 1
\}$ by the Totoki-kolmogorov  criterion, see Theorem 4.10 and
Problem 4.11 in \cite{Karatzas}. \hfill $\square$
\section{The weak solution}
\noindent Since the sequence $\{u^n,  \ n\geq 1 \}$, is tight in
$\mathcal{C}\left([0,T], \mathcal{H}\right)$, it is relatively
compact in this space. Thus, there exists a subsequence, that we
still denote by $u^n$, which converges weakly in
$\mathcal{C}\left([0,T], \mathcal{H}\right)$ to a process $u$. By
Skorohod's Representation Theorem, there exists a probability space
$\tilde{\Omega}$, and on it a sequence $\tilde{u}^n$, as well as a
noise $\tilde{F}$ and a process $\tilde{u}$ such that $(u^n , u,
F)\stackrel{\mathcal{D}}{\equiv} (\tilde{u}^n, \tilde{u} ,
\tilde{F})$  and $\tilde{u}^n$ converges almost surely to
$\tilde{u}$ in
 $\mathcal{C}\left([0,T], \mathcal{H}\right)$. We will show, by solving the
 corresponding martingale problem, that $\tilde{u}$ is a weak solution
of the equation $Eq(d, u_{0}, \sigma)$.
\begin{Pro}
For any $\varphi \in \mathcal{C}^2(D)$, such that $\varphi\equiv 0$
on $\partial D$, we have
\begin{eqnarray}
\mathcal{M}_{\varphi}(t) : & = & \int_{D}\tilde{u}(t,x)\varphi(x)dx
- \int_{D}u_{0}(x)\varphi(x)dx -
\int_{0}^{t}\int_{D}\tilde{u}(s,x)\Delta\varphi(x)dxds\nonumber\\
&& + \frac{1}{2} \sum_{i
=1}^{d}\int_{0}^{t}\int_{D}\tilde{u}^{2}(s,x)\frac{\partial
\varphi}{\partial x_i}(x)dxds  \label{chp5-5-19}
\end{eqnarray}
 is a martingale with the quadratic variation
\begin{eqnarray}\label{20}
\langle \mathcal{M}_{\varphi} \rangle_t =
\int_{0}^t\int_{D}\int_{D}\sigma(\tilde{u}(s,x))\varphi(x)f(x,y)
\sigma(\tilde{u}(s,y))\varphi(y)dxdyds.
\end{eqnarray}
\end{Pro}
\noindent\textbf{Proof}. We consider the scheme (\ref{4}) and by
 multiplying  both its sides by $\frac{1}{n^d}\varphi(x_{\kk}^n)$
 and summing over all the elements of the grid $D_n^d$, we obtain
 \begin{eqnarray*}
 \mathcal{M}_{\varphi}^n(t) :  & = & \frac{1}{n^d}\sum_{x_{\kk}^n}u^n(t,x_{\kk}^n)\varphi(x_{\kk}^n)
  -  \frac{1}{n^d}\sum_{x_{\kk}^n}u(0,x_{\kk}^n)\varphi(x_{\kk}^n)
  -\int_{0}^{t} \frac{1}{n^d} \sum_{x_{\kk}^n}\Delta^n u^n(s,x_{\kk}^n)\varphi(x_{\kk}^n)ds\\
  &  & -
  \frac{1}{2}\int_{0}^{t} \frac{1}{n^d}
  \sum_{x_{\kk}^n}\sum_{i=1}^{d}\nabla^n_i (u^n(s,x_{\kk}^n))^2\varphi(x_{\kk}^n)ds\\
 & = & \frac{1}{n^d}\sum_{x_{\kk}^n}\tilde{u}^n(t,x_{\kk}^n)\varphi(x_{\kk}^n) -
  \frac{1}{n^d}\sum_{x_{\kk}^n}u(0,x_{\kk}^n)\varphi(x_{\kk}^n)
  -\int_{0}^{t} \frac{1}{n^d} \sum_{x_{\kk}^n} \tilde{u}^n(s,x_{\kk}^n)\Delta^n\varphi(x_{\kk}^n)ds \\
  &  & + \frac{1}{2}\int_{0}^{t} \frac{1}{n^d}
  \sum_{x_{\kk}^n}\sum_{i=1}^{d}(\tilde{u}^n(s,x_{\kk}^n))^2\nabla^n_i \varphi(x_{\kk}^n)ds\\
 &  = & \int_0^t \sum_{x_{\kk}^n}\varphi(x_{\kk}^n) \sigma(\tilde{u}^n(s,x_{\kk}^n))F(ds, \II^n_{\kk}).
 \end{eqnarray*}
\noindent Note that  $\{ \mathcal{M}_{\varphi}^n(t),  t\geq 0 \}$ is
a martingale as it is a finite sum of martingales. Moreover, by
(\ref{19}) and the boundedness of $\sigma$, we can write
\begin{eqnarray*}
 \mathbb{E}(\mathcal{M}_{\varphi}^n(t))^2 & = &\sum_{x_{\kk}^n,
x_{\jj}^n}\int_{0}^{t}ds\int_{\II^n_{\kk}}dx\int_{\II^n_{\jj}} dy
\varphi(x)\sigma(\tilde{u}(s,x))\varphi(y)\sigma(\tilde{u}(s,y))f(x,y)\\
& \leq & t \sum_{x_{\kk}^n,
x_{\jj}^n}\int_{\II^n_{\kk}}dx\int_{\II^n_{\jj}} dy
\varphi(x)f(x,y)\varphi(y)\\
& \leq & t \int_{D}dx\int_{D} dy \varphi(x)f(x,y)\varphi(y).
\end{eqnarray*}
Hence, there exists a constant $c$ depending only on $ \varphi$, $k$
and $T$ such that
 $$
  \sup_{0 \leq t \leq T}\sup_{n}
 \mathbb{E}(\mathcal{M}_{\varphi}^n(t))^2 \leq  c < \infty.
 $$
 Consequently,
 $\{ \mathcal{M}_{\varphi}^n,  n\geq 1 \} $ converges  a.s.,
as $n\longrightarrow \infty$, to the martingale
  $\mathcal{M}_{\varphi}$
  given by (\ref{chp5-5-19}).
Also,  we have $$ \langle \mathcal{M}_{\varphi} \rangle _t
 = \lim_{n\longrightarrow\infty} \langle
  \mathcal{M}_{\varphi}^n\rangle_t ,
  $$
where
\begin{eqnarray*}
 \langle
  \mathcal{M}_{\varphi}^n\rangle_t
 & = & \left\langle \sum_{x_{\kk}^n}\int_0^{\cdot} \int_{\II^n_{\kk}}\varphi(x)
\sigma{(\tilde{u}^n(s,x))}F(ds,dx)\right\rangle_t\\
& = & \int_{0}^{t}\sum_{x_{\kk}^n,
x_{\jj}^n}\int_{\II^n_{\kk}}\int_{\II^n_{\jj}}\varphi(\kappa_n(x))
\sigma{(\tilde{u}^n(s,\kappa_n(x)))}f(x,y) \varphi(\kappa_n(y))
\sigma{(\tilde{u}^n(s,\kappa_n(y)))}dsdxdy,
\end{eqnarray*}
which converges a.s.  to (\ref{20}) as $n\longrightarrow \infty$.
\hfill $\square$ \newline

\noindent\textbf{Proof of the Theorem \ref{18}} Following Walsh
\cite{Walsh}, there exists a martingale measure $M$ with quadratic
variation
$$\nu(dxdyds) = \sigma{(\tilde{u}(s,x))}f(x,y)\sigma{(\tilde{u}(s,y))}dxdyds,$$
 which corresponds to the quadratic  variation
 $\langle \mathcal{M}_{\varphi} \rangle_t$. Furthermore, we
consider a measure martingale $N$ independent of $M$ and
characterized by (\ref{19}), and we set
\begin{eqnarray*}
F(t,\varphi) & = &
\int_{0}^{t}\int_{D}\frac{1}{\sigma{(\tilde{u}(s,x))}}
\mathbf1_{\{\tilde{u}(s,x)\not\in\{0,1\}\}}\varphi(x)M(dx ds)  +
\int_{0}^{t}\int_{D}\mathbf1_{\{\tilde{u}(s,x)\in\{0,1\}\}}\varphi(x)
N(dx ds).
\end{eqnarray*}
\noindent Therefore, $F(t,\varphi)$ satisfies
\begin{displaymath}
\langle F(\cdot,\varphi) \rangle_t =
\int_{0}^{t}\int_{D}\int_{D}\varphi(x)f(x,y)\varphi(y)dxdyds.
\end{displaymath}

\noindent Hence, because of the form of its quadratic variation, it
has the same distribution as $N$, and by Proposition 2.5.7 in
\cite{Anja-these}, we can write
\begin{displaymath}
 \mathcal{M}_{\varphi}(t) =
  \int_{0}^{t}\int_{D}{\sigma{(\tilde{u}(s,x))}}\varphi(x)F(dx, ds),
\end{displaymath}
\noindent therefore, $\tilde{u}$ is a weak solution of the equation
$Eq(d, u_0, \sigma)$. \hfill $\square$

\section{Appendix }
\noindent Firstly, we recall the Lemma 1.0 in \cite{Gall},
\begin{Lem}\label{3-7}
Let $Z\equiv\{Z(t), t\geq 0\}$ be a real-valued semi-martingale.
Suppose that there exists a function $\rho :
[0,+\infty)\longrightarrow[0,+\infty)$ such that
$$\int_{0}^{\varepsilon}\frac{dx}{\rho(x)} = +\infty,$$ for all
$\varepsilon >0$, and
$$\int_{0}^{t}\frac{\mathbf1_{\{Z_{s}>0\}}}{\rho(Z_{s})}d\langle
Z\rangle_{s} < +\infty,$$ for all $t>0$ a.s. Then, the local time of
 $Z$ at zero, $L_{t}^{0}(Z)$ is identically zero for all $t>0$ a.s.
\end{Lem}
\vspace{0.5cm}

\noindent The following lemma gives some useful estimates  satisfied
by $\bar{p}^n_d$.
\begin{Lem}\label{9}
\noindent i)  There exists a constant  $c>0$ such that for any $t>0$
and $\lambda>0$,
$$ \sup_{n\geq 1}\sup_{x\in
D}\|\bar{p}^n_1(t,x,\cdot)\|_{1} \leq ct^{-\lambda}\exp(-ct)
$$
\noindent ii) There exist constants $c>0$, $1 < q < 2 $ such that
for any $t>0$ and  $x, y \in D$,
$$
\sup_{n\geq 1}\int_{0}^{t}\left(\int_{D}\left|\bar{p}^n_d(t -s,x,z)
- \bar{p}^n_d(t - s, y, z)\right|dz\right)^{2}ds \leq c
\|x-y\|^{\frac{1}{q} - \frac{1}{2}};
$$
 this estimate holds with $q=1$ when $ d = 1$.\\
\noindent iii) For any $T>0$, there exists a constant $c>0$ such
that for $h>0$
\begin{eqnarray*}
\sup_{n\geq 1}\sup_{x\in
D}\sup_{t\in[0,T]}\int_{0}^{t}\|\bar{p}^n_d(t -s,x,\cdot) -
\bar{p}^n_d(t + h - s, x, \cdot)\|_{1}ds \leq c h^{\frac{1}{2}}
\end{eqnarray*}
\noindent iv) For any $T>0$, there exists a constant $c>0$ such that
for $h>0$
\begin{eqnarray*}
\sup_{n\geq 1}\sup_{x\in D}\sup_{t\in[0,T]}\int_{0}^{t}
\|\bar{p}^n_d(t-s,x,\cdot)-\bar{p}^n_d(t-s+h,x,\cdot)\|_{\alpha}^{2}ds
\leq c h^{1-\frac{\alpha}{2}}
\end{eqnarray*}
\noindent where $0 < \alpha < 2\wedge d$, and $\|\cdot\|_{\alpha}$
is defined on an appropriate class of functions $\varphi$ by
\begin{eqnarray*}
\|\varphi\|_{(\alpha)}^{2} =
\int_{D}\int_{D}|\varphi(x)||x-y|^{-\alpha}|\varphi(y)| dxdy
\end{eqnarray*}
\end{Lem}
\noindent\textbf{Proof.} The estimates $i)$, $iii)$ and $iv)$ were
proved in \cite{MLuc}, that are respectively, (A.16) of the Lemma
A.4, (A.36) and (A.38) of the Lemma A.6. Concerning $ii)$, using the
fact that $\bar{p}^n_d(t, x, z) = \prod_{i=1}^{d}\bar{p}^n_1(t -s,
x_i, z_i)$, where  $x=(x_i, 1\leq i\leq d)$ and $z = (z_i, 1\leq
i\leq d)$, the inequality $\left|\prod_{i=1}^{d}a_{i} -
\prod_{i=1}^{d}b_{i} \right| \leq \sum_{i=1}^{d} |a_{i} -
b_{i}|\prod_{j=i+1}^{d}|a_j|\prod_{j=1}^{i-1}|b_j| $ for all real
numbers $a_1, ..., a_d$, $b_1..., b_d$, the estimate $i)$ and the
Cauchy-Schwarz inequality, we can write
\begin{eqnarray*}
&& \int_{0}^{t}\left(\int_{D}\left|\bar{p}^n_d(t -s,x,z) -
\bar{p}^n_d(t - s, y, z)\right|dz\right)^{2}ds\\
&& = \int_{0}^{t}\left(\int_{D}\left|\prod_{j=1}^{d}\bar{p}^n_1(t
-s,x_j,z_j) -  \prod_{j=1}^{d}\bar{p}^n_1(t
-s,y_j,z_j)\right|dz\right)^2ds\\
&& \leq
\int_{0}^{t}\left(\sum_{i=1}^{d}\prod_{j=1}^{i-1}\|\bar{p}^n_1(t
-s,x_j,\cdot)\|_1\times \int_{0}^{1}|\bar{p}^n_1(t -s,x_i, z_i) -
\bar{p}^n_1(t -s,y_i, z_i)|dz_i\right.\\
&& \qquad\times\left.\prod_{j=i+1}^{d}\|\bar{p}^n_1(t
-s,y_j,\cdot)\|_1\right)^2ds\\
&& \leq c \int_{0}^{t}\left(\sum_{i=1}^{d}(t-s)^{-\lambda(d-1)}
\int_{0}^{1}|\bar{p}^n_1(t -s,x_i, z_i) - \bar{p}^n_1(t -s,y_i,
z_i)|dz_i \right)^2ds\\
&& \leq c \sum_{i=1}^{d} \int_{0}^{t} (t-s)^{-2\lambda(d-1)}
\int_{0}^{1}|\bar{p}^n_1(t -s,x_i, z_i) - \bar{p}^n_1(t -s,y_i,
z_i)|^2 dz_i ds.
\end{eqnarray*}
\noindent On the other hand, using notations of Gy\"{o}ngy
\cite{Gy1} and H\"{o}lder's inequality, we get
\begin{eqnarray*}
&& \int_{0}^{t} (t-s)^{-2\lambda(d-1)} \int_{0}^{1}|\bar{p}^n_1(t
-s,x_i, z_i) - \bar{p}^n_1(t -s,y_i, z_i)|^2 dz_i ds\\
&& \leq \sum_{j=1}^{n-1} \int_{0}^{t} (t-s)^{-2\lambda(d-1)}
\exp(-8j^2(t-s))ds|\varphi_j^n(x_i) - \varphi_j^n(y_i)|^2\\
&& \leq c
\sum_{j=1}^{n-1}\left(\int_{0}^{t}\exp(-8q(t-s)j^2)\right)^{1/q}(j^2|x_i
- yi|^{2}\wedge 1)\\
&& \leq c \sum_{j=1}^{n-1}\frac{1}{j^{2/q}}(j^2|x_i
- yi|^{2}\wedge 1)\\
&& \leq c |x_i - y_i|^{\frac{1}{q} - \frac{1}{2}}.
\end{eqnarray*}
\noindent Note that we choose $\lambda$ and $q$ such that $1<q<2$
and $ 0< 2\lambda (d-1)\xi <1$, where  $\frac{1}{q} + \frac{1}{\xi}
=1$. \hfill $\square$

 \vspace{1cm}
\noindent We have also used the following result in Section
\ref{sec3-3}, for the proof see Theorem 3.10 of Chapter 5 in
\cite{Ethier-Kurtz}.

\begin{Theo}\label{th4-app}
 Consider the  SDE
\begin{equation}\label{4-18}
 X_{t} = X_{0} + \int_{0}^{t}b(s, X_s)ds + \int_{0}^{t}\sigma(s, X_s)dM_s,
\end{equation}
where $t\geq 0$, $b : \R^{+}\times \R^{d} \longrightarrow\R^{d}$ and
$\sigma : \R^{+}\times\R^{d} \longrightarrow\R^{d}\otimes\R^{d}$,
are Borel continuous functions, $M$ is a $d$-dimensional martingale
in $\mathcal{M}_{c}^{2}$ and  $X$ is a $d$-dimensional process.
Assume further that the $\mathbb{R}^{d}$--valued random variable
$X_{0}$ is $\mathcal{F}_{0}$--measurable such that
$\mathbb{P}X_{0}^{-1} = \mu$. If there exists a constant $c$ such
that
\begin{equation*}
\Vert b(t,x)\Vert \leq c(1 + \Vert x \Vert) \qquad \text{and} \qquad
\left(\sum_{i,j =1}^d\sigma_{i,j}^{2}(t,x)\right)^{\frac{1}{2}}\leq
c(1 + \Vert x \Vert),
\end{equation*}
for all $t\geq 0$ and $x\in\Rd$, then for all initial condition
probability measure $\mu$ on $\Rd$, there exists a weak solution of
the SDE (\ref{4-18}).
\end{Theo}


\end{document}